\documentclass[12pt]{article} 
\usepackage{amsmath}
\usepackage{amssymb} 
\usepackage{amscd}
\usepackage{amsthm}
\usepackage[utf8]{inputenc}
\usepackage{pb-diagram}
\usepackage{hyperref}

\newtheorem{theorem}{Theorem}

\newtheorem{lemma}[theorem]{Lemma}
\newtheorem{corollary}[theorem]{Corollary}

\def\dim{{\mbox{dim}}}

\def\cala{{\mathcal A}}

\def\calr{{\mathcal R}}

\def\bbbone{\mbox{\rm 1\hspace {-.6em} l}}

\numberwithin{equation}{section}

\begin{document}

\enlargethispage{3cm}

\thispagestyle{empty}
\begin{center}
{\bf QUADRATIC DIFFERENTIAL ALGEBRAS}\\
\vspace{0.2cm}
{\bf GENERATED BY EUCLIDEAN SPACES}
\end{center}
   
\vspace{0.3cm}

\begin{center}
Michel DUBOIS-VIOLETTE
\footnote{Laboratoire de Physique Th\'eorique, UMR 8627, 
CNRS et Universit\'e Paris-Sud 11,
B\^atiment 210, F-91 405 Orsay Cedex, France\\
Michel.Dubois-Violette$@$u-psud.fr} and 
Giovanni LANDI \footnote{Matematica, Universit\`a di Trieste, Via A. Valerio, 12/1, 34127 Trieste, Italy\\ Institute for Geometry and Physics (IGAP) Trieste, Italy\\ and INFN, Trieste, Italy\\
landi$@$units.it}\end{center}
\vspace{0,5cm}


\begin{abstract}
We define a class of quadratic differential algebras which are generated as  differential graded algebras by the elements of an Euclidean space. Such a differential algebra is a differential calculus over the quadratic algebra of its elements of differential degree zero. This generalizes for arbitrary quadratic algebras the differential graded algebra of exterior polynomial differential forms for the algebra of polynomial functions on $\mathbb R^n$. We investigate the structure of these differential algebras and their connection with the Koszul complexes of quadratic algebras.
\end{abstract}

\vfill
 
 \newpage
\tableofcontents

\section{Introduction}
In this paper we shall be concerned with quadratic algebras $\cala$ generated by $n$ elements $x^\lambda$ $(\lambda\in \{1,\dots, n\})$ with relations
\[
x^\lambda x^\mu=R^{\lambda\mu}_{\nu\rho} x^\nu x^\rho
\]
for $\lambda, \mu\in \{1,\dots,n\}$ where the $R^{\lambda\mu}_{\nu\rho}$ are real numbers for $\lambda, \mu,\nu, \rho \in \{1,\dots, n\}$.\\

Let us show that for any real quadratic algebra generated by $n$ elements $x^\lambda$, one may write its relations as above. Indeed let us introduce on the real vector space $E=\sum_\lambda \mathbb R x^\lambda$ of its generators the unique Euclidean structure such that $(x^\lambda)$ is an orthonormal basis and let $\calr\subset E\otimes E$ be the space of its relations. Then $E\otimes E$ is also canonically an Euclidean space and let $P$ be the orthogonal projection onto its subspace $\calr$. One has
\[
P=\frac{1}{2}(\bbbone -R)
\]
where $R$ is an involutive symmetric endomorphism of $E\otimes E$ with matrix element $R^{\lambda\mu}_{\nu\rho}$. Then the space of relations $\calr$ is spanned by the
\[
x^\lambda x^\mu-R^{\lambda\mu}_{\nu\rho} x^\nu x^\rho
\]
with $\lambda,\mu\in \{1,\dots, n\}$ so one may write the relations of $\cala$ as expected above with furthermore the properties
\[
R^{\lambda\mu}_{\alpha \beta} R^{\alpha\beta}_{\nu\rho}= \delta^\lambda_\nu \delta^\mu_\nu,\>\>\>\> R^{\lambda\mu}_{\nu\rho}=R^{\nu\rho}_{\lambda\mu}
\]
for $\lambda, \mu, \nu, \rho\in \{1,\dots, n\}$. Since all these properties will be useful, one will start in the following by real quadratic algebras generated by an Euclidean space $E$ which is therefore not a restriction.\\

Assuming $\cala$ given as above, one can define a differential calculus over $\cala$ by setting
\begin{equation}\label{rdiff}
\left\{
\begin{array}{lll}
  x^\lambda x^\mu & =  &  R^{\lambda\mu}_{\nu\rho} x^\nu x^\rho \\
  \\
   x^\lambda dx^\mu & =  &  R^{\lambda\mu}_{\nu\rho} dx^\nu x^\rho \\
   \\
  dx^\lambda x^\mu & =  &  R^{\lambda\mu}_{\nu\rho} x^\nu dx^\rho \\
\\
   dx^\lambda dx^\mu & =  &  - R^{\lambda\mu}_{\nu\rho} dx^\nu dx^\rho \\
\end{array}
\right.
\end{equation}
which define a quadratic algebra (generated by the $x^\lambda, dx^\mu$) on which $x^\lambda \mapsto dx^\lambda$ extends uniquely as an antiderivation $d$ of square 0 $(d^2=0)$, i.e. as a differential, for the graduation given by the degree in the $dx^\lambda$. The corresponding differential graded algebra $\Omega(\cala)$  is a differential calculus over $\cala$ in the sense of \cite{wor:1989} (see also in  \cite{mdv:2001}). \\

It is worth noticing here that in the case where the $R^{\lambda\mu}_{\nu\rho}$ are given by $\delta^\lambda_\rho \delta^\mu_\nu$ then the algebra $\cala$ is the polynomial algebra $\mathbb R[x^\lambda]=S(\mathbb R^n)$ and $\Omega(\cala)$ is the algebra $S(\mathbb R^n)\otimes \wedge(\mathbb R^n)$ of polynomial differential forms on $\mathbb R^n$ while $d$ is the usual exterior differential. In the following, this case will be called {\sl the classical case}.\\

Our aim in this paper is to investigate the properties of such differential graded algebras.\\

The above construction has been introduced and used in \cite{lan-pag:2018} within a specific framework in which the symmetry condition for the $R^{\lambda\mu}_{\nu\rho}$ is not satisfied but is replaced by conditions which are adapted for a complex $\ast$-algebraic variant of the previous real formulation (see Section \ref{StarAlg}). \\

\noindent\underbar{Notations}. We use throughout the Einstein convention of summation on the repeated up-down indices. Given a finite-dimensional vector space $E$, we denote its tensor algebra by $T(E)$, its symmetric algebra by $S(E)$ and its exterior algebra by $\wedge(E)$. All the graded algebras involved in this paper are $\mathbb N$-graded algebras but it is convenient to consider these graded algebras as $\mathbb Z$-graded algebras such that their components of negative degree vanish.

\section{Canonical form of a quadratic relation}

Throughout this paper $E$ is a finite-dimensional Euclidean space with orthonormal coordinate basis $(x^\lambda)$ $\lambda\in  \{1,\dots, n\}$, $n=\dim(E)$.\\

Let $\calr$ be a linear subspace of $E\otimes E$ and let $\cala$ be the quadratic algebra defined by
\begin{equation}\label{A}
\cala=T(E)/(\calr)
\end{equation}
where $(\calr)$ denotes the two-sided ideal of the tensor algebra $T(E)$ of $E$ generated by $\calr$. Since the relations of $\cala$ are homogeneous (quadratic) of degree strictly greater than one, $\cala$ is a graded algebra
\[
\cala=\oplus_{n\in \mathbb N} \cala_n
\]
for the degree induced by the one of $T(E)$ which is unital and connected, that is 
\[
\cala_0=\mathbb R\cdot \bbbone
\]
where $\bbbone$ is the unit of $\cala$ induced from the one of $T(E)$, and one has $\cala_1=E$ while $\cala_2=E\otimes E/\calr$.\\

The space $E\otimes E$ is canonically an Euclidean space and we denote by $P$ the orthogonal projection onto the subspace $\calr$ of $E\otimes E$. One has
\begin{equation}\label{R}
P=\frac{1}{2}(\bbbone-R)
\end{equation}
where $R$ is a symmetric involution of $E\otimes E$. In terms of components in the orthonormal system $(x^\lambda\otimes x^\mu)$ one has
\begin{equation}\label{mR}
R=(R^{\lambda\mu}_{\nu\rho})
\end{equation}
with
\begin{equation}\label{invol}
R^{\lambda\mu}_{\alpha\beta}\>\> R^{\alpha\beta}_{\nu\rho}=\delta^\lambda_\nu\>\> \delta^\mu_\rho
\end{equation}
and
\begin{equation}\label{sym}
R^{\lambda\mu}_{\nu\rho}=R^{\nu\rho}_{\lambda\mu}
\end{equation}
for the real numbers $R^{\lambda\mu}_{\nu\rho}$. The quadratic algebra $\cala$ is then the real unital associative algebra generated by the $x^\lambda$ $(\lambda\in \{1,\dots, n\})$  with relations
\begin{equation}\label{cqr}
x^\lambda x^\mu - R^{\lambda\mu}_{\nu\rho}\>\> x^\nu x^\rho=0
\end{equation}
for $\lambda,\mu \in \{1,\dots,n\}$.\\

As pointed out in the introduction in the classical case, i.e. when $R^{\lambda\mu}_{\nu\rho}$ coincides with $\delta^\lambda_\rho \delta^\mu_\nu$, $\cala$ is the polynomial algebra $\mathbb R[x^\lambda]=S(\mathbb R^n)=S(E)$.

\section{The Koszul dual $\cala^!$ of $\cala$}

Let $E^\ast$ be the dual vector space of $E$ and $(\theta_\lambda)$ be the dual basis of $(x^\lambda)$. The Koszul dual $\cala^!$ of $\cala$ is the quadratic algebra  \cite{pri:1970}, \cite{man:1988}
\begin{equation}\label{Akd}
\cala^!=T(E^\ast)/(\calr^!)
\end{equation}
where $\calr^!$ is the orthogonal in $E^\ast\otimes E^\ast=(E\otimes E)^\ast$ of the subspace $\calr$ of $E\otimes E$. In terms of components, $\cala^!$ is the real unital associative algebra generated by the $\theta_\lambda$ with relations
\begin{equation}\label{cqrd}
\theta_\lambda \theta_\mu + R^{\nu\rho}_{\lambda\mu} \>\> \theta_\nu \theta_\rho=0
\end{equation}
for $\lambda,\mu\in \{1,\dots,n\}$.\\

In view of the symmetry of the $R^{\lambda\mu}_{\nu\rho}$, that is of relations (\ref{sym}), $\cala^!$ is isomorphic to the real unital algebra $\cala'$ generated by elements $y^\lambda$ for $\lambda\in\{1,\dots,n\}$ with relations
\begin{equation}\label{cqrprim}
y^\lambda y^\mu + R^{\lambda \mu}_{\nu\rho}\>\> y^\nu y^\rho=0
\end{equation}
for $\lambda,\mu\in \{1,\dots, n\}$. This latter algebra is the quadratic algebra
\begin{equation}\label{Aprim}
\cala'=T(E)/(\calr^\perp)
\end{equation}
where $\calr^\perp$ is the orthogonal subspace to $\calr$ in the Euclidean space $E\otimes E$. Note that the orthogonal projection onto  $\calr^\perp$ is 
\begin{equation}\label{-R}
P=\frac{1}{2}(\bbbone +R)
\end{equation}
thus one passes from $\cala$ to $\cala'$ by changing $R$ into $-R$ and the $x^\lambda$ into the $y^\lambda$.\\
Notice that with the identification (\ref{Aprim}),  $(y^\lambda)$ is the orthonormal basis of $E$.\\

It is not hard to see that in the classical case $\cala'$ coincides with the exterior algebra $\wedge(\mathbb R^n)=\wedge(E)$.

\section{The bigraded algebra $A$}

Let $A$ be the quadratic algebra generated by the $2n$ elements $x^\lambda,y^\mu$ with relations 
\begin{equation}\label{rA}
\left\{
\begin{array}{llll}
x^\lambda x^\mu & - & R^{\lambda\mu}_{\nu\rho}\>\> x^\nu x^\rho  &= 0   \\
\\
 x^\lambda y^\mu &  - &  R^{\lambda\mu}_{\nu\rho}\>\> y^\nu x^\rho   &  = 0 \\
 \\
 y^\lambda x^\mu &   - &  R^{\lambda\mu}_{\nu\rho}\>\> x^\nu y^\rho & = 0\\
 \\
   y^\lambda y^\mu &  + &  R^{\lambda\mu}_{\nu\rho}\>\> y^\nu y^\rho & = 0  
\end{array}
\right.
\end{equation}
for $\lambda,\mu \in \{ 1,\dots, n\}$. Since these relations are separately homogeneous in the $x^\lambda$ and in the $y^\mu$, the algebra $A$ is a bigraded algebra, that is one has
\begin{equation}\label{bigrA}
\begin{array}{l}
A = \oplus_{r,s\geq 0}\>\> A^{(r,s)}\\
\\
A^{(p,q)} A^{(r,s)} \subset A^{(p+r, q+s)}
\end{array}
\end{equation}
where in $(r,s)$, $r$ denotes the degree in the $x^\lambda$ while $s$ denotes the degree in the $y^\mu$.\\

The subalgebra of $A$ given by $\oplus_r A^{(r,0)}$ coincides with $\cala$ while the subalgebra given by $\oplus_s A^{(0,s)}$ coincides with $\cala' (\simeq \cala^!)$ that is
\begin{equation}
\begin{array}{lll}
\cala & = & \oplus_r A^{(r,0)} \subset A\\
\\
\cala' & = & \oplus_s A^{(0,s)} \subset A
\end{array}
\end{equation}
by definition.\\

\begin{lemma}\label{pprod}
In $A$, one has
\[
A^{(r,s)}=A^{(r,0)} A^{(0,s)}=A^{(0,s)}A^{(r,0)}
\]
for any $r,s\in \mathbb N$
\end{lemma}

\noindent\underbar{Proof}. By using the second and the third equalities of \ref{rdiff}, one can write any product of $(r)$ $s's$ and $(s)y's$ as a linear combination of products where the $s's$ are on the left-hand side (resp. right-hand side) of the $y's$. This implies $A^{(r,s)}\subset A^{(r,0)} A^{(0,s)}\subset A^{(r,s)}$ and $A^{(r,0)}\subset A^{(r,s)}$ and  therefore $A^{r,s)}=A^{(r,0)} A^{(0,s)}=A^{(0,s)} A^{(r,0)}.\>\square$

In other words one has
\begin{equation}\label{phr}
A=\cala \cala'=\cala'\cala
\end{equation}
in the bigraded algebra $A$ and 
\begin{equation}\label{strgr}
A^{(p,q)} A^{(r,s)}=A^{(p+r,q+s)}
\end{equation}
follows easily for any $p,q,r,s\in \mathbb N$.\\

In the classical case one has $A^{(r,s)}=\cala_r\otimes \cala'_{ s}$ since then the $x^\lambda$ and the $y^\mu$ commute, and therefore in this case $A^{(r,s)}$ coincides with $S^r(E)\otimes \wedge^s(E)=S^r(\mathbb R^n)\otimes \wedge^s(\mathbb R^n)$.

\section{The differential calculus $\Omega(\cala)$ on $\cala$}

Let us define the graded algebra $\Omega(\cala)$ by setting
\begin{equation}\label{omega}
\Omega(\cala)=\oplus_{k\in\mathbb N} \Omega^k(\cala)
\end{equation}
with
\begin{equation}\label{omegak}
\Omega^k(\cala)=\oplus_{r\in \mathbb N} A^{(r,k)}
\end{equation}
then the linear mapping $d$ of $A^{(1,0)}$ into $A^{(0,1)}$ and of $A^{(0,1)}$ into $A^{(0,2)}$ defined by
\begin{equation}\label{d}
\left\{
\begin{array}{lllll}
 x^\lambda &\mapsto   &  dx^\lambda & = & y^\lambda \\
 y^\lambda & \mapsto  &  dy^\lambda & = & 0  
\end{array}
\right.
\end{equation}
extends uniquely as an antiderivation of degree 1 of the graded algebra $\Omega(\cala)$ in view of the relations  (\ref{rA}). It follows that one has
\begin{equation}\label{diffo}
d^2=0
\end{equation}
for the derivation $d^2$ since it vanishes on the generators of $\Omega(\cala)$ (=$A$ as algebra). Thus $\Omega(\cala)$ endowed with the differential $d$ is a differential graded algebra which is generated by $E$ as graded differential algebra with the quadratic relation (\ref{rdiff}) of the introduction. Notice that the first equality of (\ref{phr}) reads
\begin{equation}\label{prd}
\Omega^k(\cala)=\omega_{\lambda_1\cdots \lambda_k}\>\> dx^{\lambda_1}\dots dx^{\lambda_k}
\end{equation}
for any $k\in \mathbb N$ with $\omega_{\lambda_1\cdots \lambda_k}\in \cala$ in view of (\ref{omegak}).\\

In the classical case the differential graded algebra $\Omega(\cala)$, that is $\Omega^\bullet(\cala)=S(\mathbb R^n)\otimes \wedge^\bullet(\mathbb R^n)$ endowed with $d$, is the differential graded algebra of exterior polynomial differential forms on $\mathbb R^n$.

\section{The codifferential $\delta$ of $\Omega(\cala)$}

By using again the relations (\ref{rA}), it follows that the linear mapping  $\delta$ of $A^{(1,0)}$ into $A^{(2,-1)}=0$ and of $A^{(0,1)}$ into $A^{(1,0)}$ defined by
\begin{equation}\label{delta}
\left\{
\begin{array}{llll}
x^\lambda  & \mapsto  &\delta x^\lambda  & =  0 \\
 y^\lambda  & \mapsto  &\delta y^\lambda  & =  x^\lambda
\end{array}
\right.
\end{equation}
extends uniquely as an antiderivation of degree -1 of the graded algebra $\Omega(\cala)$. Then since the derivation $\delta^2$ is vanishing on the generators of $\Omega(\cala)$ ($=A$ as algebra), one has
\begin{equation}\label{codiffo}
\delta^2=0
\end{equation}
that is $\delta$ is a differential of degree -1 of $\Omega(\cala)$ which will be referred to as the {\sl codifferential of} $\Omega(\cala)$.\\

Since by definition $\Omega(\cala)=A$ as algebra $\Omega(\cala)$ is bigraded and its {\sl differential degree} for which $d$ and $\delta$ are antiderivations is the second degree of this bigraduation, one has
\[
\Omega^p(\cala)=\sum_r\Omega^{(r,p)}(\cala)=\sum_r A^{(r,p)}
\]
for the differential degree $p$. For this bigraduation $d$ is of bidegree $(-1,1)$ while $\delta$ is of bidegree $(1,-1)$.\\

In the classical case the complex $(\Omega^\bullet(\cala),\delta)$, that is 
\[
\Omega^\bullet(\cala)=S(\mathbb R^n)\otimes \wedge^\bullet(\mathbb R^n)
\]
 endowed with $\delta$, identifies canonically with the Koszul complex of the quadratic algebra $\cala=S(\mathbb R^n)$.\\

\section{A generalization of the Poincar\'e lemma}

The derivation $d\delta+\delta d$ is of bidegree $(0,0)$ and one has the following lemma.
\begin{lemma}\label{todeg}
The derivation $d\delta +\delta d$ is the total degree that is one has
\begin{equation}\label{derdeg}
(d\delta +\delta d) a = (r+s) a
\end{equation}
for $a\in A^{(r,s)}$ and $r,s \in \mathbb N$.
\end{lemma}
\noindent\underbar{Proof}. This follows immediately from
\[
\left\{
\begin{array}{lll}
(d\delta + \delta d) x^\lambda  &  = & x^\lambda\\
(d\delta + \delta d) y^\mu  &  = & y^\mu
\end{array}
\right.
\]
for the generators $x^\lambda, y^\mu$, ($\lambda, \mu\in \{1,\dots, n\}$).$\>\>\square$
\begin{corollary}\label{triv}
The homology $H^\bullet(d)=\oplus_{k\in \mathbb N} H^k(d)$ of $(\Omega(\cala),d)$ is given by
\[
H^k(d)=0\>\> \text{for}\>\> k\geq 1\>\> \text{and}\>\> H^0(d)=\mathbb R
\]
while the homology $H_\bullet(\delta)=\oplus_{k\in \mathbb N}H_k(\delta)$ of $(\Omega(\cala),\delta)$ is given by 
\[
H_k(\delta)=0\>\> \text{for}\>\> k\geq 1\>\> \text{and}\>\> H_0(\delta)=\mathbb R.
\]
\end{corollary}

\noindent\underbar{Proof}. In view of (\ref{derdeg}), if $r+s>0$ then $da=0$ is equivalent to $a=d\delta \left(\frac{a}{r+s}\right)$ while $\delta a=0$ is equivalent to $a=\delta d\left(\frac{a}{r+s}\right)$. On the other hand $r+s=0$ means that $a=\lambda\bbbone$ for some $\lambda\in \mathbb R$. $\square$\\

Notice that $(\Omega(\cala),d)$ is a cochain complex while $(\Omega(\cala),\delta)$ is a chain complex, it is why the homology $H^\bullet(d)$ of $d$ is noted with the degree up while the homology $H_\bullet(\delta)$ is noted with the degree down. Both $H^\bullet(d)$ and $H_\bullet(\delta)$ are real graded algebras which here are trivial that is reduced to the ground field $\mathbb R$ (in view of above corollary).\\

Corollary (\ref{triv}) means in particular that the sequence of $\cala$-modules
\begin{equation}\label{ResR}
\dots \stackrel{\delta}{\rightarrow}\Omega^k(\cala)\stackrel{\delta}{\rightarrow} \Omega^{k-1}(\cala)\stackrel{\delta}{\rightarrow}\dots \stackrel{\delta}{\rightarrow} \Omega^1(\cala)\stackrel{\delta}{\rightarrow}\cala\stackrel{\varepsilon}{\rightarrow}\mathbb R\rightarrow 0
\end{equation}
is a resolution of the trivial module $\mathbb R$, $\varepsilon$ being the projection on degree 0 of the graded algebra $\cala$. A priori, this is not a sequence of free (i.e. projective here see in \cite{car:1958}) $\cala$-modules. However in the classical case $\cala=S(\mathbb R^n)$, (\ref{ResR})  is the Koszul resolution of $\mathbb R$ by the free $S(\mathbb R^n)$-modules $S(\mathbb R^n)\otimes \wedge^\bullet(\mathbb R^n)$.\\

It is our aim in the following to discuss this point and the relation of $(\Omega(\cala),\delta)$ with the Koszul complex of the quadratic algebra $\cala$, \cite{man:1988}.

\section{The Koszul complex $K(\cala)$}

By the very definition of the Koszul dual algebra $\cala^!$ of $\cala$, the dual vector spaces $\cala^{!\ast}_p$ of the homogeneous components $\cala^!_p$ of $\cala^!$ are given by
\begin{equation}\label{kn}
\cala^{!\ast}_p=\cap_{0\leq s\leq n-2} E^{\otimes s}\otimes \calr \otimes E^{\otimes p-s-2}
\end{equation}
for $p\geq 2$ and $\cala^{!\ast}_1=E$, $\cala^{!\ast}_0=\mathbb R\bbbone$, (as for any quadratic algebra with vector space of generators $E$ and vector space of relations $\calr \subset E\otimes E)$.\\

Notice that one has
\begin{equation}\label{aen}
\cala^{!\ast}_p\subset E^{\otimes p}
\end{equation}
for any $p\in \mathbb N$ and that furthermore
\begin{equation}\label{aepn}
\cala^{!\ast}_p\subset \calr \otimes E^{\otimes p-2}
\end{equation}
for $p\geq 2$. Therefore the linear mappings of $\cala\otimes E^{\otimes p}$ into $\cala\otimes E^{p-1}$ defined by
\[
a \otimes (x_1\otimes \dots \otimes x_p) \mapsto ax_1 \otimes (x_2\otimes \dots \otimes x_p)
\]
induces linear mappings
\begin{equation}\label{bord}
b: \cala\otimes \cala^{!\ast}_p \rightarrow \cala \otimes \cala^{!\ast}_{p-1}
\end{equation}
satisfying $b^2=0$ (with the convention that $\cala^{!\ast}_p=0$ for $p<0$). This defines the Koszul complex $K(\cala)$ of $\cala$ as 
\begin{equation}\label{K(cala)}
\dots \stackrel{b}{\rightarrow} \cala \otimes \cala^{!\ast}_p \stackrel{b}{\rightarrow} \cala\otimes \cala^{!\ast}_{p-1} \stackrel{b}{\rightarrow} \cala \otimes \cala^{!\ast}_{p-2} \stackrel{b}{\rightarrow}\dots
\end{equation}
which ends as
\begin{equation}\label{K0}
\dots \stackrel{b}{\rightarrow} \cala \otimes \calr \stackrel{b}{\rightarrow} \cala \otimes E \stackrel{b}\rightarrow \cala
\end{equation}
since $\cala^{!\ast}_2 = \calr$ and $\cala^{!\ast}_1=E$. As well known and easy to show the exactness of the following sequence
\begin{equation}\label{pres}
\cala\otimes \calr \stackrel {b}{\rightarrow} \cala \otimes E\stackrel{b}{\rightarrow}\cala \stackrel{\varepsilon}{\rightarrow} \mathbb R \rightarrow 0
\end{equation}
is equivalent to the presentation of $\cala$ by generators and relations. Thus if the Koszul complex $K(\cala)$ is acyclic in positive degrees, (i.e. in degrees $\geq 2$ since by (\ref{pres}) it is acyclic in degree =1), it gives a free resolution
\begin{equation}\label{kreso}
K(\cala) \stackrel{\varepsilon}{\rightarrow} \mathbb R \rightarrow 0
\end{equation}
of the trivial module $\mathbb R$. A quadratic algebra $\cala$ is said to be a {\sl Koszul algebra}  \cite{pri:1970},  \cite{man:1988} whenever its Koszul complex $K(\cala)$ is acyclic in positive degrees and then the resolution (\ref{kreso}) is refered to as the {\sl Koszul resolution} of $\mathbb R$. In this case, the Koszul resolution is a minimal free resolution.

\section{Comparison of $K(\cala)$ with $(\Omega(\cala),\delta)$}

There are surjective homomorphisms
\begin{equation}\label{hp}
h_p:\cala\otimes \cala^{!\ast}_p\rightarrow \Omega^p(\cala)=\cala \cdot \cala'_p
\end{equation}
of (left) $\cala$-modules which are induced by the isomorphisms $\cala^{!\ast}_p\simeq \cala'_p$ of vector spaces. However the corresponding diagram of homomorphisms of $\cala$-modules
\begin{equation}\label{DKO}
\begin{diagram}
\node{\dots} \arrow{e,t}{b} \node{K_p(\cala)}\arrow{s,r}{h_p} \arrow{e,t}{b}\node{K_{p-1}(\cala)}\arrow{s,r}{h_{p-1}}
\arrow{e,t}{b} \node{K_{p-2}(\cala)}\arrow{s,r}{h_{p-2}}
\arrow{e,t}{b}\node{\dots}\\
\node{\dots} \arrow{e,t}{\delta} \node{\Omega^p(\cala)}\arrow{e,t}{\delta}\node{\Omega^{p-1}(\cala)}\arrow{e,t}{\delta} \node{\Omega^{p-2}(\cala)}\arrow{e,t}{\delta}\node{\dots}
\end{diagram}
\end{equation}
is generally not commutative and, in contrast to the $K_p(\cala)=\cala\otimes \cala^{!\ast}_p$ the $\Omega^p(\cala)=\cala \cdot \cala'_p$ are generally not free modules.\\

However it is worth noticing here that the exact sequence (\ref{pres}) corresponding to the presentation of $\cala$ coincides canonically with the end
\[
\Omega^2(\cala) \stackrel{\delta}{\rightarrow} \Omega^1(\cala) \stackrel{\delta}{\rightarrow} \cala \stackrel{\varepsilon}{\rightarrow} \mathbb R \rightarrow 0
\]
of the exact sequence (\ref{ResR}) since one has $\Omega^1(\cala)=\cala\otimes E$ and $\Omega^2(\cala)= \cala\otimes \cala'_2\simeq \cala\otimes \calr$ canonically in view of $\cala'_2=E\otimes E/\calr^\perp$ and $E\otimes E=\calr\oplus \calr^\perp$ and since furthermore with these identifications $\delta$ coincide with $b$ up to an irrelevant factor. Thus the end of (\ref{DKO})
\[
\begin{diagram} 
\node{K_2(\cala)}\arrow{s,r}{h_2} \arrow{e,t}{b} \node{K_1(\cala)}\arrow{s,r}{h_1} \arrow{e,t}{b}\node{\cala}\arrow{s,r}{h_0}
\arrow{e,t}{\varepsilon } \node{\mathbb R}\arrow{s}\arrow{s,r}{=}
\arrow{e}\node{0}\\
\node{\Omega^2(\cala)} \arrow{e,t}{\delta} \node{\Omega^1(\cala)}\arrow{e,t}{\delta}\node{\cala}\arrow{e,t}{\varepsilon} \node{\mathbb R}\arrow{e}\node{0}
\end{diagram}
\]
is a commutative diagram of $\cala$-modules with vertical isomorphisms. It is also worth noticing that in the classical case where $\cala =S(\mathbb R^n)$ and $R^{\lambda\mu}_{\nu\rho}=\delta^\lambda_\rho \delta^\mu_\nu$, the diagram (\ref{DKO}) is commutative with vertical isomorphisms. It is therefore natural to ask for conditions on $(\cala, R)$ for which this occurs. In view of the acyclicity of (\ref{ResR}), this can occur only if $\cala$ is Koszul and if the $\Omega^p(\cala)$ are free modules. A very natural condition under which this is satisfied is when $R:E\otimes E \rightarrow E\otimes E$ satisfies the so called quantum Yang-Baxter equation
\begin{equation}\label{YB}
(R\otimes I)(I\otimes R)(R\otimes I)=(I\otimes R)(R\otimes I)(I\otimes R)
\end{equation}
as endomorphim of $E\otimes E \otimes E$, where $I$ denotes the identity mapping of $E=\mathbb R^n$ into itself. Indeed it is well known that then $\cala$ is a Koszul algebra \cite{gur:1990}, \cite{wam:1993} and that by using the fact that the $I^{\otimes^r} \otimes R\otimes I^{\otimes^s}$ generate a representation of the symmetric group  $\frak S(r+s+2)$ in $E^{\otimes r+s+  2}$ one can show that the product $\cala \cdot \cala'$ in $A$ is isomorphic to $\cala\otimes \cala'$ and that the $b's$ coincides then to the $\delta's$ up to unimportant normalizations.

\section{The complex quadratic $\ast$-algebra variant}\label{StarAlg}

The variant described in this section is useful and natural for a quantum theoretical formalism and for discussion of reality conditions in a noncommutative geometrical framework. It consists here in considering complex quadratic algebras which are complex $\ast$-algebras generated by hermitian elements which belong to an Euclidean subspace generating a Hilbertian space.\\

So let again $E$ be a $n$-dimensional Euclidean space with orthonormal coordinate basis $(x^\lambda)\>\> \lambda\in \{1,\dots, n\}$ and let $E_c$ be the complexified of $E$ endowed with the corresponding sesquilinear product. Thus $E_c$ is a Hilbertian space endowed with an antilinear involution $x\mapsto x^\ast$ with $E=\{x\in E_c\vert x^\ast=x\}$ and $(x^\lambda)$ is an orthonormal basis of $E_c$ which consists of hermitian elements. Then there is a unique extension of the involution of $E_c$ to the complex tensor algebra $T(E_c)$ for which $T(E_c)$ is a complex $\ast$-algebra.\\

Let $\calr$ be a $\ast$-invariant subspace of $E_c\otimes E_c$, (i.e. invariant by $x\otimes y\mapsto y^\ast\otimes x^\ast)$, the complex quadratic algebra
\[
\cala=T(E_c)/(\calr)
\]
is a $\ast$-algebra which is generated by the hermitian elements $x^\lambda,\>\> \lambda\in \{1,\dots,n\}$. The space $E_c\otimes E_c$ is canonically a Hilbertian space and let $P$ be the hermitian projector onto $\calr\subset E_c \otimes E_c$. One has
\[
P=\frac{1}{2}(\bbbone -R)
\]
where $R$ is a hermitian involution. In terms of components in the $x^\lambda\otimes x^\mu$ one has $R=(R^{\lambda\mu}_{\nu\rho})$ with $R^{\lambda\mu}_{\nu\rho}\in \mathbb C$ satisfying 
\[
R^{\lambda\mu}_{\alpha\beta} R^{\alpha\beta}_{\nu\rho} = \delta^\lambda_\nu \delta^\mu_\rho
\]
 as before and
\begin{equation}\label{astR}
R^{\lambda\mu}_{\nu\rho} = \bar R^{\nu\rho}_{\lambda\mu} = \bar R^{\mu\lambda}_{\rho\nu}
\end{equation}
instead of relation (\ref{sym}) of Section 2. The relations of $\cala$ reads 
\[
x^\lambda x^\mu= R^{\lambda\mu}_{\nu\rho} x^\nu x^\rho
\]
for $\lambda, \mu \in \{1,\dots, n\}$ and by the same construction as in Sections 4, 5, 6, 7 using a complex $R$ as above instead of the real one there, one constructs a differential calculus $\Omega(\cala)$ for $\cala$ which is a quadratic complex differential graded $\ast$-algebra generated by $E$ and one defines the codifferential $\delta$ on $\Omega(\cala)$, etc. 
This variant of the construction has been used in a particular 4-dimensional setting in \cite{lan-pag:2018}.

\section{Generalizations}

Finally it is worth noticing that Relations (\ref{rdiff}) define a quadratic differential graded algebra for any $R=(R^{\lambda\mu}_{\nu\rho})$ which is involutive $(R^2=\bbbone)$ or more generally which is invertible and satisfies an equation of the form
\begin{equation}\label{pHe}
R^2=\alpha R + \beta \bbbone
\end{equation}
where $\alpha$ and $\beta$ are two scalar (with $\beta\not = 0$), in order that the equations with $R$ replaced by $R^n$ with $n\geq 2$ do not imply other constraints.\\

Then all the above analysis applies, even the one of Section 9 with some slight modifications in view of the work of \cite{gur:1990} and \cite{wam:1993} implying the Koszulity of $\cala$ under the assumption (\ref{YB}), (Hecke case).

\section*{Acknowlegements}
This work was partially supported by INDAM-GNSAGA and by the INDAM-CNRS Project LIA-LYSM.

\newpage

\end{document}